\documentclass[reqno,12pt]{amsart}
\usepackage{amsmath, epsfig, cite}
\usepackage{amssymb}
\usepackage{amsfonts}
\usepackage{latexsym}
\usepackage{amsthm}

\newtheorem{theorem}{Theorem}%[section]

\newtheorem{conjecture}{Conjecture}
\newtheorem{lemma}{Lemma}

\theoremstyle{remark}

\numberwithin{equation}{section}

\newcommand*{\pFq}[5]{{}_{#1}F_{#2}\!\left[ \begin{matrix} #3\\[5pt] #4\end{matrix};#5\right]}

\allowdisplaybreaks[1]

\author{Victor J.\ W.\ Guo}
\address{School of Mathematics and Statistics, Huaiyin Normal University,
Huai'an 223300, Jiangsu, People's Republic of China}
\email{jwguo@hytc.edu.cn}
\thanks{The first author was partially supported by the National Natural
Science Foundation of China (grant 11771175).}

\author{Ji-Cai Liu}
\address{Department of Mathematics, Wenzhou University, Wenzhou 325035,
People's Republic of China}
\email{jcliu2016@gmail.com}
\thanks{The second author was partially supported by the National Natural
Science Foundation of China (grant 11801417).}

\author{Michael J.\ Schlosser}
\address{Fakult\"at f\"ur Mathematik, Universit\"at Wien,
Oskar-Morgenstern-Platz~1, A-1090 Vienna, Austria}
\email{michael.schlosser@univie.ac.at}
\thanks{The third author was partially supported by FWF Austrian Science
Fund grant P 32305.}
\title{An extension of a supercongruence of Long and Ramakrishna}



\begin{document}

\begin{abstract}
We prove two supercongruences for specific truncated hypergeometric series.
These include an uniparametric extension of a supercongruence that was recently
established by Long and Ramakrishna.
Our proofs involve special instances of various hypergeometric identities including
Whipple's transformation and the Karlsson--Minton summation.
\end{abstract}

\maketitle

\section{Introduction}
Let $(a)_n=a(a+1)\cdots(a+n-1)$ denote the Pochhammer symbol.
For complex numbers $a_0,a_1,\ldots,a_r$ and $b_1,\ldots,b_r$,
the (generalized) hypergeometric series $_{r+1}F_r$ is defined as
$$
{}_{r+1}F_{r}\!\left[\begin{array}{cccc}
a_0, & a_1, & \ldots,     & a_{r}\\[5pt]
     & b_1, & \ldots,     & b_{r}
\end{array};z
\right]
=\sum_{k=0}^\infty\frac{(a_0)_k (a_1)_k\cdots (a_r)_k}
{k!(b_1)_k\cdots (b_r)_k}z^k.
$$
Summation and transformation formulas for generalized hypergeometric series
play an important part in the investigation of supercongruences.
See, for instance, \cite{GL2,Liu0,Liu,LR,MO,Mortenson4,Wang}. In particular,
Long and Ramakrishna~\cite[Theorems~3 and 2]{LR}
proved the following two supercongruences:
\begin{equation}
\sum_{k=0}^{p-1}\frac{\left(\frac{1}{2}\right)_k^3}{k!^3}
\equiv
\begin{cases} -\displaystyle \Gamma_p\!\left(\frac{1}{4}\right)^4
\pmod{p^3}, &\text{if $p\equiv 1\pmod 4$,}\\[10pt]
 -\frac{p^2}{16}\displaystyle
 \Gamma_p\!\left(\frac{1}{4}\right)^4\pmod{p^3},
&\text{if $p\equiv 3\pmod 4$,}
\end{cases} \label{eq:lr3}
\end{equation}
and
\begin{equation}
\sum_{k=0}^{p-1} (6k+1) \frac{\left(\frac{1}{3}\right)_k^6}{k!^6}
\equiv
\begin{cases} -p\,\displaystyle \Gamma_p\!\left(\frac{1}{3}\right)^9
\pmod{p^6}, &\text{if $p\equiv 1\pmod 6$,}\\[10pt]
 -\frac{10}{27}p^4\,\displaystyle
 \Gamma_p\!\left(\frac{1}{3}\right)^9\pmod{p^6},
&\text{if $p\equiv 5\pmod 6$,}
\end{cases} \label{eq:d2}
\end{equation}
where $\Gamma_p(x)$ is the $p$-adic Gamma function.
The restriction of the supercongruence \eqref{eq:lr3} modulo $p^2$
was earlier established by Van Hamme \cite[Equation (H.2)]{Hamme}.
The supercongruence \eqref{eq:d2} is even stronger than
a conjecure made by Van Hamme \cite[Equation (D.2)]{Hamme}
who asserted the corresponding supercongruence modulo
$p^4$ for $p\equiv 1\pmod 6$.
Long and Ramakrishna also mentioned that \eqref{eq:d2} does not hold
modulo $p^7$ in general.

The first purpose of this paper is to prove the following supercongruence.
Note that the $r=\pm 1$ cases partially confirm the $d=5$ and $q\to 1$ case
of \cite[Conjectures~1 and 2]{GS19}.

\begin{theorem}\label{thm:1}
  Let $r\leqslant 1$ be an odd integer coprime with $5$.
  Let $p$ be a prime such that
$p\equiv-\frac{r}{2}\pmod{5}$ and $p\geqslant \frac{5-r}{2}$. Then
\begin{equation}
  \sum_{k=0}^{p-1}(10k+r)\frac{\left(\frac{r}{5}\right)_k^5}{k!^5}
  \equiv 0 \pmod{p^4},
\label{eq:thm1}
\end{equation}
\end{theorem}

Recently, the second author \cite{Liu-2020}
established the following supercongruence related to \eqref{eq:d2}:
\begin{equation}
\sum_{k=0}^{p-1}(6k-1)\frac{\left(-\frac{1}{3}\right)_k^6}{k!^6}\equiv
\begin{cases}
  \displaystyle140p^4\,\Gamma_p\!\left(\frac{2}{3}\right)^9\pmod{p^5},
  \quad &\text{if $p\equiv 1\pmod{6}$,}\\[15pt]
  \displaystyle378p\,\Gamma_p\!\left(\frac{2}{3}\right)^{9}\pmod{p^5},
  \quad &\text{if $p\equiv 5\pmod{6}$,}
\end{cases}\label{a-1}
\end{equation}
where $p$ is a prime.

The second purpose of this paper is to give the following common generalization
of the second supercongruence in \eqref{eq:d2}, restricted to modulo $p^5$,
and the first supercongruence in \eqref{a-1}.

\begin{theorem}\label{thm:2}
  Let $r\leqslant 1$ be an integer coprime with $3$. Let $p$ be
  a prime such that $p\equiv-r\pmod{3}$ and $p\geqslant 3-r$. Then
\begin{align}
\sum_{k=0}^{p-1} (6k+r) \frac{\left(\frac{r}{3}\right)_k^6}{k!^6}
&\equiv \frac{(-1)^{r+1}\,80rp^4}{81}\cdot
\frac{\Gamma_p\!\left(1+\frac{r}{3}\right)^2}
{\Gamma_p\!\left(1+\frac{2r}{3}\right)^3
\Gamma_p\!\left(1-\frac{r}{3}\right)^4}\notag\\
  &\quad\;\times\sum_{k=0}^{1-r}\frac{(r-1)_k\left(\frac{r}{3}\right)_k^3}
    {(1)_k\left(\frac{2r}{3}\right)_k^3} \pmod{p^5}. \label{eq:thm2}
\end{align}
\end{theorem}
Letting $r=1$ and $r=-1$ in \eqref{eq:thm2} and using \eqref{bb-2} and
\eqref{bb-4}, we arrive at the $p\equiv 5\pmod{6}$ case of \eqref{eq:d2}
modulo $p^5$ and the $p\equiv 1\pmod{6}$ case of \eqref{a-1}, respectively.

Our proof of Theorem \ref{thm:1} will require Whipple's well-poised
$_7F_6$ transformation formula (see \cite[p.~28]{Bailey}):
\begin{align}
&{}_7F_{6}\!\left[\begin{array}{ccccccc}
a, & 1+\frac{1}{2}a, & b,     & c,      & d,      & e,     & -n \\[5pt]
   & \frac{1}{2}a,   & 1+a-b, & 1+a-c,  & 1+a-d,  & 1+a-e, & 1+a+n
\end{array};1
\right]  \notag\\[5pt]
&=\frac{(a+1)_n (a-d-e+1)_n}{(a-d+1)_n (a-e+1)_n}\,
{}_4F_{3}\!\left[\begin{array}{cccc}
1+a-b-c,       &  d,       &  e,        & -n     \\[5pt]
 d+e-a-n,     &  1+a-b,   & 1+a-c  &
\end{array};1
\right],  \label{eq:6f5-0}
\end{align}
where $n$ is a non-negative integer, and Karlsson--Minton's summation
formula (see, for example, \cite[Equation~(1.9.2)]{GR}):
\begin{equation}
{}_{r+1}F_{r}\!\left[\begin{array}{cccc}
-n, & b_1+m_1, & \ldots,     & b_{r}+m_r\\[5pt]
     & b_1, & \ldots,     & b_{r}
\end{array};1
\right]
=0,  \label{eq:km}
\end{equation}
where $n,m_1,\ldots,m_r$ are non-negative integers and $n>m_1+\cdots+m_r$.
Our proof of Theorem~\ref{thm:2} relies on a $_7F_6$ transformation formula
slightly different from Whipple's $_7F_6$ transformation formula
\eqref{eq:6f5-0}, obtained as a result from combining \eqref{eq:6f5-0}
with a $_4F_3$ transformation formula. The transformation was already utilized by the second author to prove
\eqref{a-1}.

Furthermore, in order to prove Theorem \ref{thm:2},
we require some properties of the $p$-adic Gamma function,
collected in the following two lemmas.
\begin{lemma}\label{lem-1}{\rm\cite[Section 11.6]{cohen-b-2007}}
Suppose $p$ is an odd prime and $x\in \mathbb{Z}_p$. Then
\begin{align}
&\Gamma_p(0)=1,\quad\Gamma_p(1)=-1,\label{bb-1}\\
&\Gamma_p(x)\Gamma_p(1-x)=(-1)^{a_p(x)},\label{bb-2}\\
  &\Gamma_p(x)\equiv \Gamma_p(y)\pmod{p}
    \quad\text{for $x\equiv y\pmod{p}$},\label{bb-3}\\[10pt]
&\frac{\Gamma_p(x+1)}{\Gamma_p(x)}=
\begin{cases}
-x\quad&\text{if $v_p(x)=0$,}\\
-1\quad &\text{if $v_p(x)>0$, }
\end{cases}\label{bb-4}
\end{align}
where $a_p(x)\in \{1,2,\dots,p\}$ with $x\equiv a_p(x)\pmod{p}$ and
$v_p(\cdot)$ denotes the $p$-order.
\end{lemma}

\begin{lemma} \label{lem-2} {\rm\cite[Lemma 17, (4)]{LR}}
Let $p$ be an odd prime. If $a\in \mathbb{Z}_p, n\in \mathbb{N}$
such that none of $a,a+1,\dots,a+n-1$ are in $p\mathbb{Z}_p$, then
\begin{equation}
(a)_n=(-1)^n\frac{\Gamma_p(a+n)}{\Gamma_p(a)}.\label{bb-6}
\end{equation}
\end{lemma}

In the following Sections~\ref{sec:thm1} and
\ref{sec:thm2}, we give proofs of
Theorems~\ref{thm:1} and \ref{thm:2}, respectively.
The final Section~\ref{sec:disc} is devoted to a discussion
and includes two conjectures.

%%%%%%%%%%%%%%%%%%%%%%%%%%%%%%%%%%%%%%%%%%%%%%%%%%%%%%%%%%%%%%%%%%%%%%%%%%%%%%%%%%%%%%%%%%%%%%%%%%%%%%%%%%%%%%%%%%%%%%%%%%%%%%%%%%%%%%%%%%%%%%%%%%%%%%%%%%%%%%%%%%%%%%%%%%%%%%%%%%%%%%%%%%%%%%%%%%%%%%%%%%%%%%%%%%%%%%%%%%%%%%%%%%%%%

\section{Proof of Theorem \ref{thm:1}}\label{sec:thm1}
Motivated by the work of McCarthy and Osburn~\cite{MO} and
Mortenson~\cite{Mortenson4}, we take the following choice of parameters
in \eqref{eq:6f5-0}.
Let $a=\frac{r}{5}$, $b=\frac{r+5}{10}$, $c=\frac{r+3p}{5}$,
$d=\frac{r+3ip}{5}$, $e=\frac{r-3ip}{5}$, and $n=\frac{3p-r}{5}$,
where $i^2=-1$. Then we conclude that
\begin{align}
&{}_6F_{5}\!\left[\begin{array}{ccccccc}
     \frac{r}{5}, & 1+\frac{r}{10}, & \frac{r+3p}{5},     & \frac{r+3ip}{5},
                  & \frac{r-3ip}{5},    & \frac{r-3p}{5} \\[5pt]
                    & \frac{r}{10},   & 1-\frac{3p}{5},     & 1-\frac{3ip}{5},
                    & 1+\frac{3ip}{5},     & 1+\frac{3p}{5}
\end{array};1
\right]  \notag\\[5pt]
&=\frac{\left(1+\frac{r}{5}\right)_{\frac{3p-r}{5}}
\left(1-\frac{r}{5}\right)_{\frac{3p-r}{5}}}
{\left(1-\frac{3ip}{5}\right)_{\frac{3p-r}{5}}
\left(1+\frac{3ip}{5}\right)_{\frac{3p-r}{5}} }\,
{}_4F_{3}\!\left[\begin{array}{cccc}
                   \frac{5-r-6p}{10}, & \frac{r+3ip}{5}, & \frac{r-3ip}{5},
                   & \frac{r-3p}{5} \\[5pt]
\frac{2r-3p}{5},        &  \frac{r+5}{10},       & \frac{5-3p}{5}
\end{array};1
\right],  \label{eq:6f5}
\end{align}

It is easy to see that, for $k\geqslant 0$ and any $p$-adic integer $b$,
\begin{equation}
  \left(a+bp\right)_k \left(a-bp\right)_k \left(a+bip\right)_k
  \left(a-bip\right)_k \equiv (a)_k^4 \pmod{p^4}.  \label{eq:modp4}
\end{equation}
Hence, the left-hand side of \eqref{eq:6f5} is congruent to
\begin{align*}
\sum_{k=0}^{\frac{3p-r}{5}}\frac{\left(1+\frac{r}{10}\right)_k
\left(\frac{r}{5}\right)_k^5}{\left(\frac{r}{10}\right)_k (1)_k^5}
&=\frac{1}{r}\sum_{k=0}^{\frac{3p-r}{5}}(10k+r)
\frac{\left(\frac{r}{5}\right)_k^5}{ k!^5} \\[5pt]
  &\equiv\frac{1}{r}\sum_{k=0}^{p-1}(10k+r)
    \frac{\left(\frac{r}{5}\right)_k^5}{ k!^5} \pmod{p^4},
\end{align*}
where we have used the fact that
$\frac{\left(\frac{r}{5}\right)_k}{k!}\equiv 0\pmod p$ for
$\frac{3p-r}{5}<k\leqslant p-1$ (the condition $p\geqslant \frac{5-r}{2}$
in the theorem is to guarantee $\frac{3p-r}{5}\leqslant p-1$).
Since $\frac{3p-r}{5}\geqslant \frac{2p+r}{5}$, we have
\begin{align*}
\frac{\left(1+\frac{r}{5}\right)_{\frac{3p-r}{5}}
\left(1-\frac{r}{5}\right)_{\frac{3p-r}{5}}}
{\left(1-\frac{3ip}{5}\right)_{\frac{3p-r}{5}}
\left(1+\frac{3ip}{5}\right)_{\frac{3p-r}{5}} }
&=\frac{\left(1+\frac{r}{5}\right)_{\frac{3p-r}{5}}
\left(1-\frac{r}{5}\right)_{\frac{3p-r}{5}}}
{\left(1+\frac{9p^2}{25}\right)_{\frac{3p-r}{5}} } \\[5pt]
&\equiv 0\pmod{p^2}.
\end{align*}
Finally, by the congruences
\begin{equation}
\left(a+bip\right)_k \left(a-bip\right)_k
\equiv  \left(a+bp\right)_k \left(a-bp\right)_k
\equiv (a)_k^2
\pmod{p^2}\label{eq:modp2}
\end{equation}
for any $p$-adic integer $b$, we obtain
\begin{align*}
{}_4F_{3}\!\left[\begin{array}{cccc}
                   \frac{5-r-6p}{10}, & \frac{r+3ip}{5}, & \frac{r-3ip}{5},
                   & \frac{r-3p}{5}    \\[5pt]
\frac{2r-3p}{5},        &  \frac{r+5}{10},       & \frac{5-3p}{5}
\end{array};1
\right]
&\equiv
{}_4F_{3}\!\left[\begin{array}{cccc}
                   \frac{5-r-6p}{10}, & \frac{r}{5}, & \frac{r}{5},
                   & \frac{r-3p}{5}    \\[5pt]
\frac{2r-3p}{5},        &  \frac{r+5}{10},       & \frac{5-3p}{5}
\end{array};1
\right]  \\[5pt]
&\equiv
{}_4F_{3}\!\left[\begin{array}{cccc}
                   \frac{5-r-6p}{10}, & \frac{r+p}{5}, & \frac{r-p}{5},
                   & \frac{r-3p}{5}    \\[5pt]
\frac{2r-3p}{5},        &  \frac{r+5}{10},       & \frac{5-3p}{5}
\end{array};1
\right] \\[5pt]
&=0 \pmod{p^2},
\end{align*}
where we have utilized Karlsson--Minton's summation \eqref{eq:km} with
$n=\frac{3p-r}{5}$, $b_1=\frac{2r-3p}{5}$, $b_2=\frac{r+5}{10}$,
$b_3=\frac{5-3p}{5}$, $m_1=\frac{1-r}{2}$, $m_2=\frac{2p+r-5}{10}$,
and $m_3=\frac{2p+r-5}{5}$ in the last step.

%%%%%%%%%%%%%%%%%%%%%%%%%%%%%%%%%%%%%%%%%%%%%%%%%%%%%%%%%%%%%%%%%%%%%%%%%%%%%%%%%%%%%%%%%%%%%%%%%%%%%%%%%%%%%%%%%%%%%%%%%%%%%%%%%%%%%%%%%%%%%%%%%%%%%%%%%%%%%%%%

\section{Proof of Theorem \ref{thm:2}}\label{sec:thm2}
We can verify \eqref{eq:thm2} for $r=1$ and $p=2$ by hand. In what follows,
we assume that $p$ is an odd prime. Recall the following transformation
formula \cite[Equation~(4.2)]{Liu-2020}:
\begin{align}
  &\pFq{7}{6}{t,&1+\frac{t}{2},&-n,&t-a,&t-b,&t-c,&1-t-m+n+a+b+c}
  {&\frac{t}{2},&1+t+n,&1+a,&1+b,&1+c,&2t+m-n-a-b-c}{1}\notag\\[10pt]
  &=\frac{(1+t)_n(a+b+2-m-t)_n(a+c+2-m-t)_n(b+c+2-m-t)_n}
    {(1+a)_n(1+b)_n(1+c)_n(a+b+c+1-m-2t)_n}\notag\\[10pt]
  &\quad\;\times\frac{(a+b+1-m-t)(a+c+1-m-t)(b+c+1-m-t)}
    {(a+b+n+1-m-t)(a+c+n+1-m-t)(b+c+n+1-m-t)}\notag\\[10pt]
  &\quad\;\times\pFq{4}{3}{-m,&\hskip -9mm -n,&\hskip -9mm a+b+c+1-m-2t,
 &a+b+c+1+n-m-t}{&\hskip -1cm a+b+1-m-t,&a+c+1-m-t,&b+c+1-m-t}{1}.\label{d-1}
\end{align}

Let $\zeta$ be a fifth primitive root of unity.
Setting $m=1-r$, $t=\frac{r}{3}$, $n=\frac{2p-r}{3}$, $a=\frac{2p\zeta}{3}$,
$b=\frac{2p\zeta^2}{3}$ and $c=\frac{2p\zeta^3}{3}$ in \eqref{d-1} and using
$1+\zeta+\zeta^2+\zeta^3+\zeta^4=0$, the left-hand side of \eqref{d-1} becomes
\begin{align*}
&\pFq{7}{6}{1+\frac{r}{6},&\frac{r}{3},&\frac{r-2p}{3},&\frac{r-2p\zeta}{3},
&\frac{r-2p\zeta^2}{3},&\frac{r-2p\zeta^3}{3},&\frac{r-2p\zeta^4}{3}}
  {&\frac{r}{6},&1+\frac{2p}{3},&1+\frac{2p\zeta}{3},&1+\frac{2p\zeta^2}{3},
  &1+\frac{2p\zeta^3}{3},&1+\frac{2p\zeta^4}{3}}{1}\\[10pt]
&\equiv \frac{1}{r}\sum_{k=0}^{\frac{2p-r}{3}}
(6k+r) \frac{\left(\frac{r}{3}\right)_k^6}{k!^6}\pmod{p^5},
\end{align*}
where we have used the facts that none of the denominators in $_7F_6$
contain a multiple of $p$ (the condition $p\geqslant 3-r$ in the theorem
is to guarantee $\frac{2p-r}{3}\leqslant p-1$) and
\begin{align*}
\left(u+vp\right)_k\left(u+vp\zeta\right)_k\left(u+vp\zeta^2\right)_k
\left(u+vp\zeta^3\right)_k\left(u+vp\zeta^4\right)_k\equiv (u)_k^5\pmod{p^5}.
\end{align*}
Furthermore, for $\frac{2p-r}{3}<k\leqslant p-1$ we have
$\left(\frac{r}{3}\right)_k \equiv 0\pmod{p}$. Thus,
\begin{align}
&\pFq{7}{6}{1+\frac{r}{6},&\frac{r}{3},&\frac{r-2p}{3},&\frac{r-2p\zeta}{3},
&\frac{r-2p\zeta^2}{3},&\frac{r-2p\zeta^3}{3},&\frac{r-2p\zeta^4}{3}}
  {&\frac{r}{6},&1+\frac{2p}{3},&1+\frac{2p\zeta}{3},&1+\frac{2p\zeta^2}{3},
  &1+\frac{2p\zeta^3}{3},&1+\frac{2p\zeta^4}{3}}{1}\notag\\[10pt]
  &\equiv \frac{1}{r}\sum_{k=0}^{p-1} (6k+r)
    \frac{\left(\frac{r}{3}\right)_k^6}{k!^6}\pmod{p^5}.\label{new-1}
\end{align}

On the other hand, we determine the terminating hypergeometric series
on the right-hand side of \eqref{d-1} modulo $p$:
\begin{align}
  &\frac{(a+b+1-m-t)(a+c+1-m-t)(b+c+1-m-t)}
    {(a+b+n+1-m-t)(a+c+n+1-m-t)(b+c+n+1-m-t)}\notag\\[10pt]
  &\times\pFq{4}{3}{-m,&\hskip -5mm -n,&\hskip -5mm a+b+c+1-m-2t,
 &a+b+c+1+n-m-t}{&\hskip -1cm a+b+1-m-t,&a+c+1-m-t,&b+c+1-m-t}{1}\notag\\[10pt]
&\equiv 8\sum_{k=0}^{1-r}\frac{(r-1)_k\left(\frac{r}{3}\right)_k^3}
{(1)_k\left(\frac{2r}{3}\right)_k^3} \pmod{p}.\label{new-2}
\end{align}
Moreover,
\begin{align}
  &\frac{(1+t)_n(a+b+2-m-t)_n(a+c+2-m-t)_n(b+c+2-m-t)_n}
    {(1+a)_n(1+b)_n(1+c)_n(a+b+c+1-m-2t)_n}\notag\\[10pt]
  &=\frac{\left(1+\frac{r}{3}\right)_{\frac{2p-r}{3}}
    \left(1+\frac{2r+2p(\zeta+\zeta^2)}{3}\right)_{\frac{2p-r}{3}}
    \left(1+\frac{2r+2p(\zeta+\zeta^3)}{3}\right)_{\frac{2p-r}{3}}
    \left(1+\frac{2r+2p(\zeta^2+\zeta^3)}{3}\right)_{\frac{2p-r}{3}}}
{(-1)^{\frac{2p-r}{3}}\left(1+\frac{2p\zeta}{3}\right)_{\frac{2p-r}{3}}
\left(1+\frac{2p\zeta^2}{3}\right)_{\frac{2p-r}{3}}
\left(1+\frac{2p\zeta^3}{3}\right)_{\frac{2p-r}{3}}
\left(1+\frac{2p\zeta^4}{3}\right)_{\frac{2p-r}{3}}}.\label{new-3}
\end{align}
Note that
\begin{equation}
  \left(1+\frac{r}{3}\right)_{\frac{2p-r}{3}}=\frac{2p}{3}
  \left(1+\frac{r}{3}\right)_{\frac{2p-r-3}{3}},\label{new-4}
\end{equation}
and
\begin{align}
&\left(1+\frac{2r+2p(\zeta+\zeta^2)}{3}\right)_{\frac{2p-r}{3}}
     \left(1+\frac{2r+2p(\zeta+\zeta^3)}{3}\right)_{\frac{2p-r}{3}}
 \left(1+\frac{2r+2p(\zeta^2+\zeta^3)}{3}\right)_{\frac{2p-r}{3}}\notag\\[10pt]
&=\frac{5p^3}{27}\left(1+\frac{2r+2p(\zeta+\zeta^2)}{3}\right)_{\frac{p-2r-3}{3}}
\left(1+\frac{2r+2p(\zeta+\zeta^3)}{3}\right)_{\frac{p-2r-3}{3}}\notag\\[10pt]
 &\quad\;\times \left(1+\frac{2r+2p(\zeta^2+\zeta^3)}{3}\right)_{\frac{p-2r-3}{3}}
    \left(\frac{3+p(2\zeta+2\zeta^2+1)}{3}\right)_{\frac{p+r}{3}} \notag\\[10pt]
  &\quad\;\times \left(\frac{3+p(2\zeta+2\zeta^3+1)}{3}\right)_{\frac{p+r}{3}}
    \left(\frac{3+p(2\zeta^2+2\zeta^3+1)}{3}\right)_{\frac{p+r}{3}}.\label{new-5}
\end{align}
Combining \eqref{new-4} and \eqref{new-5}, we arrive at
\begin{align}
&\frac{\left(1+\frac{r}{3}\right)_{\frac{2p-r}{3}}
\left(1+\frac{2r+2p(\zeta+\zeta^2)}{3}\right)_{\frac{2p-r}{3}}
\left(1+\frac{2r+2p(\zeta+\zeta^3)}{3}\right)_{\frac{2p-r}{3}}
\left(1+\frac{2r+2p(\zeta^2+\zeta^3)}{3}\right)_{\frac{2p-r}{3}}}
{(-1)^{\frac{2p-r}{3}}\left(1+\frac{2p\zeta}{3}\right)_{\frac{2p-r}{3}}
\left(1+\frac{2p\zeta^2}{3}\right)_{\frac{2p-r}{3}}
\left(1+\frac{2p\zeta^3}{3}\right)_{\frac{2p-r}{3}}
\left(1+\frac{2p\zeta^4}{3}\right)_{\frac{2p-r}{3}}}\notag\\[10pt]
&\equiv \frac{(-1)^{\frac{2p-r}{3}}10p^4}{81}\cdot
\frac{\left(1+\frac{r}{3}\right)_{\frac{2p-r-3}{3}}
\left(1+\frac{2r}{3}\right)^3_{\frac{p-2r-3}{3}}
(1)^3_{\frac{p+r}{3}}}{(1)^4_{\frac{2p-r}{3}}}\pmod{p^5}.\label{new-6}
\end{align}
It follows from \eqref{new-1}--\eqref{new-3} and \eqref{new-6} that
\begin{align}
\sum_{k=0}^{p-1} (6k+r) \frac{\left(\frac{r}{3}\right)_k^6}{k!^6}
&\equiv \frac{(-1)^{\frac{2p-r}{3}}\,80rp^4}{81}\cdot
\frac{\left(1+\frac{r}{3}\right)_{\frac{2p-r-3}{3}}
\left(1+\frac{2r}{3}\right)^3_{\frac{p-2r-3}{3}}(1)^3_{\frac{p+r}{3}}}
{(1)^4_{\frac{2p-r}{3}}}\notag\\
&\quad\;\times\sum_{k=0}^{1-r}\frac{(r-1)_k\left(\frac{r}{3}\right)_k^3}
{(1)_k\left(\frac{2r}{3}\right)_k^3}\pmod{p^5}.
\label{new-7}
\end{align}

By Lemmas \ref{lem-1} and \ref{lem-2}, we have
\begin{align}
&\frac{\left(1+\frac{r}{3}\right)_{\frac{2p-r-3}{3}}
\left(1+\frac{2r}{3}\right)^3_{\frac{p-2r-3}{3}}
(1)^3_{\frac{p+r}{3}}}{(1)^4_{\frac{2p-r}{3}}}\notag\\[10pt]
&\overset{\eqref{bb-6}}{=}\frac{(-1)^{\frac{2p-r}{3}+r}\,
\Gamma_p\!\left(\frac{2p}{3}\right)
\Gamma_p\!\left(\frac{p}{3}\right)^3
\Gamma_p\!\left(1+\frac{p+r}{3}\right)^3\Gamma_p(1)}
{\Gamma_p\!\left(1+\frac{r}{3}\right)
\Gamma_p\!\left(1+\frac{2r}{3}\right)^3
\Gamma_p\!\left(1+\frac{2p-r}{3}\right)^4}\notag\\[10pt]
&\overset{\eqref{bb-3}}{\equiv}
\frac{(-1)^{\frac{2p-r}{3}+r}\,
\Gamma_p\!\left(0\right)^4\Gamma_p\left(1+\frac{r}{3}\right)^2\Gamma_p(1)}
{\Gamma_p\!\left(1+\frac{2r}{3}\right)^3
\Gamma_p\!\left(1-\frac{r}{3}\right)^4}\pmod{p}\notag\\[10pt]
&\overset{\eqref{bb-1}}{=}
\frac{(-1)^{\frac{2p-r}{3}+r+1}\,
\Gamma_p\!\left(1+\frac{r}{3}\right)^2}
 {\Gamma_p\!\left(1+\frac{2r}{3}\right)^3
 \Gamma_p\!\left(1-\frac{r}{3}\right)^4}.\label{new-8}
\end{align}
The proof of \eqref{eq:thm2} then follows from \eqref{new-7} and \eqref{new-8}.

%%%%%%%%%%%%%%%%%%%%%%%%%%%%%%%%%%%%%%%%%%%%%%%%%%%%%%%%%%%%%%%%%%%%%%%%%%%%%%%%%%%%%%%%%%%%%%%%%%%%%%%%%%%%%%%%%%%%%%%%%%%%%%%%%%%%%%%%%%%%%%%%%%%%%%%%%%%%%%%%

\section{Discussion}\label{sec:disc}

We know that many supercongruences have nice $q$-analogues
(see \cite{Guo-rima,Guo-a2,GS19,GS1,GuoZu,LP,Zu19}).
For example, we have the following conjectural $q$-analogue of \eqref{eq:thm1}:
for the same $p$ and $r$ as in Theorem \ref{thm:1},
\begin{equation}
  \sum_{k=0}^{p-1}[10k+r]\frac{\left(q^r;q^5\right)_k^5}
  {\left(q^5;q^5\right)_k^5}q^{\frac{5(3-r)k}{2}}
\equiv 0 \pmod{[p]^4},
\label{eq:conj}
\end{equation}
where $[n]=1+q+\cdots+q^{n-1}$ is the $q$-integer and
$(a;q)_n=(1-a)(1-aq)\cdots (1-aq^{n-1})$ denotes the $q$-shifted factorial.

Although there are $q$-analogues of Whipple's well-poised $_7F_6$
transformation and of Karlsson--Minton's summation
(see \cite[Appendix (II.27) and (III.18)]{GR}),
we are unable to give a proof of \eqref{eq:conj}.
This is because we only know a $q$-analogue of \eqref{eq:modp2}
(see \cite[Lemma~1]{GS1}) but do not know any $q$-analogues of
\eqref{eq:modp4}. Besides, we do not know how to prove \eqref{eq:conj} by
using the method of `creative microscoping' devised in \cite{GuoZu} either.

While in Theorem~\ref{thm:2} we were able to provide a common generalization
of the second supercongruence  in \eqref{eq:d2} (restricted to modulo $p^5$)
and the first supercongruence in \eqref{a-1}, it appears to be rather
difficult to extend Theorem~\ref{thm:1} to a higher supercongruence
involving the $p$-adic Gamma function in the spirit of Theorem \ref{thm:2}, even in the special cases $r=1$ or $r=-1$.

%While in Theorem~\ref{thm:2} we were able to provide a common generalization
%of the second supercongruence  in \eqref{a-0}
%(restricted to modulo $p^5$) and the first supercongruence in \eqref{a-1},
%it appears to be very difficult to give a similar result that would
%be a common generalization of the first supercongruence in \eqref{a-0}
%(restricted to modulo $p^5$)
%and the second supercongruence in \eqref{a-1}.

We end our paper with two further conjectures for future research.
Conjecture~\ref{conj-one} concerns a stronger version of
Theorem~\ref{thm:2} and includes the second supercongruence in
\eqref{eq:d2} as a special case.
Conjecture~\ref{conj-three} concerns a common generalization of the
first supercongruence in \eqref{eq:d2} and the second
supercongruence in \eqref{a-1}.

\begin{conjecture}\label{conj-one}
The supercongruence \eqref{eq:thm2} holds modulo $p^6$ for any prime $p>3$.
\end{conjecture}

\begin{conjecture}\label{conj-three}
Let $r\leqslant 1$ be an integer coprime with $3$. Let $p\geqslant 7$
be a prime such that $p\equiv r\pmod{3}$ and $p\geqslant 3-2r$. Then
\begin{align}
\sum_{k=0}^{p-1} (6k+r) \frac{\left(\frac{r}{3}\right)_k^6}{k!^6}
&\equiv \frac{(-1)^r\,8rp}{3}\cdot
\frac{\Gamma_p\!\left(1+\frac{r}{3}\right)^2}
{\Gamma_p\!\left(1+\frac{2r}{3}\right)^3
\Gamma_p\!\left(1-\frac{r}{3}\right)^4}\notag\\
  &\quad\;\times\sum_{k=0}^{1-r}\frac{(r-1)_k\left(\frac{r}{3}\right)_k^3}
    {(1)_k\left(\frac{2r}{3}\right)_k^3} \pmod{p^6}. \label{conj-3}
\end{align}
\end{conjecture}

We remark that by using \eqref{d-1} and the same method as in the proof of
Theorem~\ref{thm:2}, we can only show that \eqref{conj-3} holds modulo $p^2$.
A new technique is needed to prove Conjecture~\ref{conj-three}.

\end{document}